\documentclass[12pt]{article}
\usepackage[cp1251]{inputenc}
\usepackage[T2A]{fontenc}
\usepackage[english,ukrainian]{babel}
\usepackage{amsmath,amsfonts,amssymb,amsthm}
\usepackage[space]{cite}

\usepackage{geometry}
\usepackage[bookmarksnumbered, colorlinks, plainpages]{hyperref}

\numberwithin{equation}{section}
\newtheorem{theorem}{Теорема}[section]

\newtheorem{corollary}[theorem]{Наслідок}

\theoremstyle{remark}

\theoremstyle{definition}
\newtheorem{example}[theorem]{Приклад}
\newtheorem{definition}[theorem]{Означення}

\begin{document}

\begin{flushleft}

\medskip
\textbf{V. O. Soldatov }\small{(Inst. Math. Nat. Acad. Sci. Ukraine, Kyiv; Inst. Appl. Math. and Mech. Nat. Acad. Sci. Ukraine, Sloviansk)} \normalsize
\medskip

\large\textbf{Solvability of linear boundary-value problems for ordinary differential systems in the space $\mathbf{C^{n}}$}

\end{flushleft}

\begin{flushleft}
	
	\textbf{В. О. Солдатов }\small{(Ін-т математики НАН України, Київ; Ін-т прикл. математики і механіки НАН України, Слов'янськ)}
	
	\normalsize

	\medskip
	
	\large
	\textbf{Розв'язність лінійних крайових задач для звичайних диференціальних систем у просторі $\mathbf{C^{n}}$}
\end{flushleft}

\normalsize

\medskip

\begin{abstract}
\noindent We study linear boundary-value problems for systems of first-order ordinary differential equations with the most general boundary conditions in the normed spaces of continuously differentiable functions on a finite closed interval. The boundary conditions are allowed to be overdetermined or underdetermined with respect to the differential system and may contain arbitrary derivatives of the unknown functions. We prove that the problem operator is Fredholm on appropriate pairs of normed spaces, find its  index and $d$-characteristics, and prove limit theorems for sequences of the characteristic matrices of the boundary-value problems under study and $d$-characteristics of the corresponding Fredholm operators.

\end{abstract}


\section{Вступ}\label{section0}

Лінійні крайові задачі для систем звичайних диференціальних рівнянь широко використовують для моделювання різноманітних процесів як в природі, так і в економіці чи суспільстві. Питання щодо розв'язності таких задач та властивостей їх розв'язків є істотно більш складними, ніж для добре вивченої задачі Коші, з огляду на велику різноманітність крайових умов.
Умови коректної розв'язності неоднорідних лінійних крайових задач із досить загальними крайовими умовами та неперервності їх розв'язків за параметром встановлено в працях І.~Т.~Кігурадзе \cite{Kiguradze1975, Kiguradze1988} і його послідовників \cite{Ashordia1996, MikhailetsReva2008DAN9, KodliukMikhailetsReva2013UMJ,  MikhailetsChekhanova2015JMathSci, MikhailetsPelekhataReva2018UMJ, PelekhataReva2019UMJ}. Характерною особливістю таких задач є задання крайових умов у вигляді $By=q$ з використанням довільного лінійного неперервного оператора $B$ на банаховому (дійсному або комплексному) просторі $r-1$ разів неперервно диференційовних функцій, де $r$~--- порядок диференціальних рівнянь, які утворюють систему.

Згодом клас досліджуваних крайових задач істотно розширили. Було уведено і досліджено лінійні крайові задачі, розв'язки яких пробігають вибраний комплексний банахів функціональний простір~--- простір $n$ разів неперервно диференційовних функцій \cite{Soldatov2015UMJ, MikhailetsMurachSoldatov2016MFAT}, де $n\geq r$, простір Гельдера порядку $s>r$ \cite{MikhailetsMurachSoldatov2016EJQTDE, MasliukSoldatov2018MFAT, Masliuk2017UMJ} або простір Соболєва--Слободецького порядку $s\geq r$ \cite{MikhailetsReva2008DAN8, KodlyukMikhailets2013JMS, GnypKodlyukMikhailets2015UMJ, Gnyp2016UMJ, GnypMikhailetsMurach2017EJDE, MasliukMikhailets2018UMJ, AtlasiukMikhailets2019UMJ11, MikhailetsSkorobohach2021UMJ}, за природного припущення, що крайовий оператор $B$ є неперервним на вибраному просторі. Для цих задач знайдено достатні та необхідні умови коректної розв'язності та неперервності за параметром розв'язків у вказаних просторах, та показано, що похибка і нев'язка розв'язку мають однаковий порядок малості.

Отримані результати знайшли застосування до дослідження граничних властивостей матриць Гріна крайових задач \cite{KodliukMikhailetsReva2013UMJ, MikhailetsChekhanova2015JMathSci}, різних багатоточкових крайових задач \cite{MikhailetsMurachSoldatov2016MFAT, Atlasiuk2020UMJ8} та їх апроксимативних властивостей щодо загальних крайових задач \cite{MasliukPelekhataSoldatov2020MFAT, MurachPelekhataSoldatov2021UMJ}.

В останній час до досліджуваного класу було долучено лінійні диференціальні системи з недоозначеними або переозначеними крайовими умовами, розв'язки яких пробігають вибраний простір Соболєва--Слободецького \cite{AtlasiukMikhailets2019UMJ10, MikhailetsAtlasiukSkorobohach2023UMJ, MikhailetsAtlasiuk2024BJMA, MikhailetsAtlasiuk2024CMP}. Було розроблено  новий підхід до аналізу таких крайових задач за допомогою їх прямокутних числових характеристичних матриць. Він дозволив встановити фредгольмовість оператора задачі, знайти його індекс і виразити $d$-характеристики оператора (тобто вимірності його ядра і коядра) у термінах характеристичної матриці задачі та, як наслідок, встановити достатню і необхідну умову коректної розв'язності задачі. Крім того, було доведено граничні теореми для послідовностей характеристичних матриць досліджуваних задач, з яких випливає напівнеперервність зверху $d$-характеристик операторів цих задач.

Мета цієї статті~--- розвинути зазначений підхід, поширивши його на лінійні диференціальні системи, розв'язки яких розглядаються у вибраному просторі неперервно диференційовних функцій. У статті розглядаються лінійні задачі з найбільш загальними крайовими умовами для систем звичайних диференціальних рівнянь першого порядку, розв'язки яких пробігають комплексний банахів простір $n$ разів неперервно диференційовних функцій на скінченному замкненому інтервалі дійсної осі, де $n$~--- довільне натуральне число. Кількість скалярних крайових умов є довільною, тобто їх набір може бути як недоозначеним, так і переозначеним щодо диференціальної системи. У статті доведено фредгольмовість оператора задачі, знайдено його індекс і $d$-характеристики, а також доведено граничні теореми для послідовностей характеристичних матриць і $d$-характеристик задач. Частину результатів проілюстровано на прикладах одноточкової і багатоточкової крайових задач.

\section{Постановка задачі}\label{section1}
Нехай числа $a,b\in\mathbb{R}$ задовольняють умову $a<b$. Позначимо коротко через $C^{(n)}$, де $n\in\mathbb{N}$, комплексний банахів простір усіх $n$ разів неперервно диференційовних на $[a,b]$ функцій $x:[a,b]\to\mathbb{C}$, наділений нормою
$$
\|x\|_{(n)}:=\sum_{j=0}^{n}\,\max\limits_{a\leq t\leq b}|x^{(j)}(t)|.
$$
Зокрема, $C=C^{(0)}$~--- банахів простір усіх неперервних на $[a,b]$ комплекснозначних функцій. Позначимо відповідно через $(C^{(n)})^{m}$ та $(C^{(n)})^{m\times m}$, де $m\in\mathbb{N}$, комплексні банахові простори, утворені вектор-функціями вимірності $m$ та матрицями-функціями розміру $m\times m$, усі компоненти яких належать до
$C^{(n)}$; при цьому вектори інтерпретуємо як стовпці. Норма вектор-функції у просторі $(C^{(n)})^{m}$ дорівнює сумі норм усіх її компонент у $C^{(n)}$, а норма матриці-функції у просторі $(C^{(n)})^{m\times m}$ дорівнює максимуму норм  усіх її стовпців у $(C^{(n)})^{m}$. Норми у цих просторах позначаємо також через $\|\cdot\|_{(n)}$, оскільки з контексту завжди буде зрозуміло про норму в якому просторі (скалярних функцій, вектор-функцій чи матриць-функцій) йде мова. Як звичайно, $\mathbb{C}^{l}$ та $\mathbb{C}^{l\times m}$~--- банахові простори усіх комплексних числових векторів вимірності $l$ та комплексних числових матриць розміру $l\times m$ відповідно, де $l,m\in\mathbb{N}$. Для числових векторів і матриць використовуємо аналогічні норми, які позначаємо через $\|\cdot\|$.

Нехай довільно вибрано числа $m,l,n\in\mathbb{N}$. На скінченному інтервалі $(a,b)$ розглянемо таку лінійну крайову задачу для системи $m$ диференціальних рівнянь першого порядку:
\begin{gather}\label{bound_pr_1}
(Ly)(t):=y^{\prime}(t) + A(t)y(t)=f(t), \quad t\in(a,b),\\
By=c. \label{bound_pr_2}
\end{gather}
Тут довільно задано матрицю-функцію $A\in(C^{(n-1)})^{m\times m}$, вектор-функцію $f\in(C^{(n-1)})^{m}$, числовий вектор $c\in\mathbb{C}^{l}$ і лінійний неперервний оператор
$$
B\colon(C^{(n)})^{m} \rightarrow\mathbb{C}^{l}.
$$
Крайова умова \eqref{bound_pr_2} задає $l$ скалярних крайових умов для системи $m$ диференціальних рівнянь першого порядку. У випадку $l>m$ ця крайова умова є \textit{перевизначеною}, а у випадку  $l<m$ вона є \textit{недовизначеною} щодо диференціальної системи \eqref{bound_pr_1}.
Розв'язок розглянутої крайової задачі інтерпретуємо як  вектор-функцію $y\in(C^{(n)})^m$, що задовольняє рівняння \eqref{bound_pr_1} на $(a,b)$ і умову \eqref{bound_pr_2}. Відображення $y\mapsto Ly$ є лінійним неперервним оператором на парі просторів $(C^{(n)})^m$ і $(C^{(n-1)})^m$.

Звісно, розв'язки системи \eqref{bound_pr_1} заповнюють увесь простір $(C^{(n)})^m$, якщо його права частина $f(\cdot)$ пробігає увесь простір $(C^{(n-1)})^m$. Отже, крайова умова \eqref{bound_pr_2} є найбільш загальною для цієї системи. Вона охоплює всі відомі типи класичних крайових умов (початкові умови задачі Коші, дво- та багатоточкові умови, інтегральні та інтегро-диференціальні умови), а також некласичні крайові умови, які містять похідні порядку $j\leq n$ шуканої функції.

\section{Основні результати}\label{section2}

Сформулюємо основні результати статті.

Перепишемо крайову задачу \eqref{bound_pr_1}, \eqref{bound_pr_2} у вигляді рівняння $(L,B)y=(f,c)$ за допомогою лінійного неперервного оператора
\begin{equation}\label{th2-LB}
(L,B)\colon (C^{(n)})^m\rightarrow (C^{(n-1)})^m\times\mathbb{C}^{l}.
\end{equation}

\begin{theorem}\label{th_fredh-bis}
Оператор \eqref{th2-LB} є фредгольмовим з індексом $m-l$.
\end{theorem}

Нагадаємо, що лінійний неперервний оператор $T\colon X \rightarrow Y$
на парі банахових просторів $X$ і $Y$ називають фредгольмовим, якщо його ядро $\ker T:=\{x\in X:Tx=0\}$ та коядро $\operatorname{coker}T:=Y/T(X)$ скінченновимірні. Якщо оператор $T$  фредгольмів, то його область визначення $T(X)$ замкнена в $Y$ \cite[Лема~19.1.1]{Hermander2007}. Індексом фредгольмового оператора $T$ називають число
$$
\mathrm{ind}\,T:=\dim\ker T-\dim\operatorname{coker}T\in\mathbb{Z}.
$$

У випадку $l=m$ теорема~\ref{th_fredh-bis} встановлена в \cite[теорема~1]{MikhailetsChekhanova2014DAN7}.

Позначимо через $Y$ єдиний розв'язок класу $(C^{(n)})^{m \times m}$ такої матричної задачі Коші:
\begin{equation}\label{Cachy}
Y'(t)+A(t)Y(t)=O_m, \quad t \in (a,b), \quad\quad Y(a)=I_m,
\end{equation}
де $O_m$ і $I_m$ позначають відповідно нульову і одиничну матриці розміру $m\times m$.

\begin{definition}\label{matrix_BY}
Комплексну числову матрицю розміру $l \times m$ називаємо \emph{характеристичною} для крайової задачі \eqref{bound_pr_1}, \eqref{bound_pr_2}, якщо кожний стовпчик цієї матриці є результатом дії оператора $B$ на стовпчик з тим же номером матриці-функції $Y$. Цю характеристичну матрицю позначаємо через $M(L,B)$.
\end{definition}

\begin{theorem}\label{th_FR_n}
Для $d$-характеристик фредгольмового оператора \eqref{th2-LB} правильні формули
\begin{gather}\label{dimker}
\dim\ker(L,B)=m-\operatorname{rank}M(L,B),\\
\dim\operatorname{coker}(L,B)=l-\operatorname{rank}M(L,B). \label{dimcoker1}
\end{gather}
\end{theorem}

Як звичайно, $\operatorname{rank}M(L,B)$ позначає ранг числової матриці $M(L,B)$. Отже, вимірності ядра і коядра оператора \eqref{th2-LB} відповідно збігаються з вимірностями ядра і коядра лінійного оператора, породженого матрицею $M(L,B)$ на парі просторів $\mathbb{C}^{m}$ і $\mathbb{C}^{l}$.

\begin{corollary}\label{th_invertible-bis}
Оператор \eqref{th2-LB} є оборотним (тобто топологічним ізоморфізмом) тоді й лише тоді, коли  $l=m$ і квадратна матриця $M(L,B)$ є невиродженою.
\end{corollary}

Поряд із задачею \eqref{bound_pr_1}, \eqref{bound_pr_2}, розглянемо послідовність крайових задач вигляду
\begin{gather}\label{Lk-sq}
(L_k y)(t):=y^{\prime}(t) + A_k(t)y(t)=f_k(t), \quad t\in(a,b),\\
\label{Bk-sq}
B_k y=c_k,
\end{gather}
параметризованих числом $k\in\mathbb{N}$. Тут $A_k\in(C^{(n-1)})^{m\times m}$, $f_k\in(C^{(n-1)})^{m}$, $c_k\in\mathbb{C}^{l}$, а $B_k$~--- лінійний неперервний оператор на парі просторів $(C^{(n)})^{m}$ і $\mathbb{C}^{l}$.

Пов'яжемо з крайовими задачами \eqref{Lk-sq}, \eqref{Bk-sq} послідовність лінійних неперервних операторів
\begin{equation}\label{LBk}
(L_k,B_k)\colon (C^{(n)})^m\rightarrow (C^{(n-1)})^m\times\mathbb{C}^l
\end{equation}
і послідовність їх характеристичних матриць $M(L_k,B_k)\in \mathbb{C}^{l\times m}$, параметризовані числом~$k$. Як звичайно, $(L_k,B_k)\xrightarrow{s}(L, B)$ позначає збіжність послідовності  операторів $(L_k,B_k)$ до $(L, B)$ у сильній операторній топології.

\begin{theorem}\label{th_hm_stc}
Якщо $(L_k,B_k)\xrightarrow{s}(L, B)$, то $M(L_k,B_k)\to M(L,B)$ (границі при $k\to\infty$).
\end{theorem}

\begin{theorem}\label{th_stc-col}
Припустимо, що
\begin{equation}\label{LB-st-conv-con}
(L_k,B_k)\xrightarrow{s}(L,B)\quad\mbox{при}\quad k\to\infty.
\end{equation}
Тоді
\begin{gather}\label{inq ker}
\operatorname{dim}\operatorname{ker}(L_k,B_k)\leq
\operatorname{dim}\operatorname{ker}(L,B)
\quad\mbox{для всіх}\quad k\gg1,\\
\operatorname{dim}\operatorname{coker}(L_k,B_k)\leq
\operatorname{dim}\operatorname{coker}(L,B)
\quad\mbox{для всіх}\quad k\gg1.
\label{inq coker}
\end{gather}
\end{theorem}

У наступних трьох наслідках останньої теореми припустимо, що виконується умова~\eqref{LB-st-conv-con}.

\begin{corollary}\label{th_stc-col-1}
Якщо оператор \eqref{th2-LB} є оборотним, то оператори \eqref{LBk} є оборотними для всіх $k\gg1$.
\end{corollary}

\begin{corollary}\label{th_stc-col-2}
Якщо крайова задача  \eqref{bound_pr_1}, \eqref{bound_pr_2} має розв'язок для кожних зазначених правих частин, то це правильно і для крайових задач \eqref{Lk-sq}, \eqref{Bk-sq} для всіх $k\gg1$.
\end{corollary}

\begin{corollary}\label{th_stc-col-3}
Якщо крайова задача \eqref{bound_pr_1}, \eqref{bound_pr_2}  має лише тривіальний розв'язок у випадку нульових правих частин, то це правильно і для крайових задач \eqref{Lk-sq}, \eqref{Bk-sq} для всіх $k\gg1$.
\end{corollary}

Доведення сформульованих теорем дамо у п.~\ref{section-proofs}. Наведені наслідки негайно випливають із відповідних теорем.

\section{Приклади}\label{section-exmpls}

\begin{example}
Розглянемо крайову задачу, яка складається з диференціальної системи \eqref{bound_pr_1} і одноточкової крайової умови
\begin{equation}\label{1.3t1}
By:=\sum_{j=0}^{n}\alpha_{j}\,y^{(j)}(a)=c,
\end{equation}
де $A\in\mathbb{C}^{m \times m}$ (тобто усі коефіцієнти диференціальної системи сталі) і кожне $\alpha_{j}\in\mathbb{C}^{l\times m}$. Відображення $y\mapsto(Ly,By)$ є лінійним неперервним оператором на парі просторів \eqref{th2-LB}. Відповідна матрична задача Коші \eqref{Cachy} має розв'язок $Y(t)=e^{(a-t)A}$. Отже,
\begin{equation*}
M(L,B)=\sum_{j=0}^{n}\alpha_{j}Y^{(j)}(a)=\sum_{j=0}^{n}\alpha_{j}(-A)^j,
\end{equation*}
що дозволяє обчислити $d$-характеристики оператора $(L,B)$ за  теоремою~\ref{th_FR_n} і зробити висновок, що вони не залежать від довжини інтервалу $(a,b)$, на якому задано цю крайову задачу (порівняти з \cite[Приклад~3.1]{MikhailetsAtlasiuk2024BJMA}).
\end{example}

\begin{example}
Розглянемо крайову задачу, яка складається з диференціальної системи \eqref{bound_pr_1} з коефіцієнтом $A(t)\equiv O_{m}$ і багатоточкової крайової умови вигляду
\begin{equation*}
By:=\sum_{j=1}^{N}\sum_{i=0}^{s_{j}} 
\beta_{j,i}\,({}^{\mathcal{C}}D_{b-}^{\alpha_{j,i}} y)(t_j)=q,
\end{equation*}
 заданої для $N\in\mathbb{N}$ точок $t_{1},\ldots,t_{N}\in [a,b]$. Тут усі  $s_{j}\in\mathbb{N}$, $\beta_{j,i}\in\mathbb{C}^{l\times m}$, а порядки $\alpha_{j,i}\in\mathbb{C}$ дробових похідних задовольняють умови $\operatorname{Re}\alpha_{j,i}\geq0$, $[\operatorname{Re}\alpha_{j,i}]\leq n-1$ (як звичайно, $[\sigma]$~--- ціла частина дійсного числа $\sigma$), а також
\begin{equation}\label{condition=0}
\alpha_{j,i}=0\;\Leftrightarrow\;i=0.
\end{equation}
Тут кожна дробова похідна ${{}^{\mathcal{C}}D}_{b-}^{\alpha_{j,i}}$ розуміється у сенсі Капуто. Неперервність оператора $B$ на парі просторів $(C^{(n)})^{m}$ і $\mathbb{C}^{l}$ випливає безпосередньо з відповідних властивостей дробових похідних Капуто і вказаних обмежень на їх порядки (див., наприклад, \cite[c.~93]{KilbasSrivastavaTrujillo2006}). Оскільки одинична матриця
$Y(\cdot)=I_m$ є розв'язком задачі Коші \eqref{Cachy}, де $A(t)\equiv O_{m}$, то 
\begin{equation*}
M(L,B)= \sum_{j=1}^{N}\sum_{i=0}^{s_{j}}\beta_{j,i}({}^{\mathcal{C}}D_{b-}^{\alpha_{j,i}} I_m)=\sum_{j=1}^{N}\beta_{j,0}
\end{equation*}
на підставі умови \eqref{condition=0}. Це дозволяє обчислити $d$-характеристики оператора $(L,B)$ за  теоремою~\ref{th_FR_n}  і зробити висновок, що вони не залежать від довжини інтервалу $(a,b)$, точок $t_{1},\ldots,t_{N}$ і матриць $\beta_{j,i}$, для яких $i\geq1$
(порівняти з \cite[Приклад~3.2]{MikhailetsAtlasiuk2024BJMA}). 
\end{example}

\section{Доведення}\label{section-proofs}

\textit{Доведення теореми $\ref{th_fredh-bis}$}. Запишемо
\begin{equation*}
(L,B)=(L,C_{l,m})+(0,B-C_{l,m}),
\end{equation*}
де $C_{l,m}$~--- довільний лінійний неперервний оператор на парі просторів $(C^{(n)})^m$ і $\mathbb{C}^{l}$. Оскільки тут другий доданок є скінченновимірним (а отже, компактним) оператором, то або оператори $(L,B)$ і $(L,C_{l,m})$ фредгольмові на парі просторів \eqref{th2-LB} і мають однаковий індекс, або вони не є фредгольмовими на цій парі  \cite[розд.~IV, теорема~5.26]{Kato1995}. Розглянемо окремо три можливі випадки: $l=m$, $l>m$ і $l<m$.

У випадку $l=m$ покладемо $C_{m,m}y:=y(a)$ для довільного $y\in(C^{(n)})^m$. Із теореми про існування та єдиність розв'язку задачі Коші для диференціальної системи \eqref{bound_pr_1} негайно випливає, що відображення $y\mapsto(Ly,C_{m,m}y)$ є ізоморфізмом на парі просторів
\begin{equation}\label{(L,C)isomorphism}
(L,C_{m,m}):(C^{(n)})^m\leftrightarrow(C^{(n-1)})^m\times\mathbb{C}^{m},
\end{equation}
а він є окремим випадком фредгольмовості оператора з нульовим індексом. Отже, оператор \eqref{th2-LB} успадковує останню властивість.

У випадку $l>m$ покладемо $C_{l,m}y:=\operatorname{col}(y(a),\theta_{l-m})$
для довільного $y\in(C^{(n)})^m$, де $\theta_{l-m}$ позначає нуль-вектор у просторі $\mathbb{C}^{l-m}$. Лінійний неперервний оператор
\begin{equation}\label{L-Cauchy}
(L,C_{l,m}):(C^{(n)})^m\to(C^{(n-1)})^m\times\mathbb{C}^{l}
\end{equation}
має нульове ядро та область значень $(C^{(n-1)})^m\times\mathbb{C}^{m}\times\{0\}^{l-m}$ на підставі щойно згаданої теореми. Отже, цей оператор разом із $(L,B)$ є фредгольмовим з індексом $0-(l-m)$.

В останньому випадку $l<m$ покладемо $C_{l,m}y:=(y_{1}(a),\ldots,y_{l}(a))$ для  довільного $y=(y_{1},\ldots,y_{m})\in(C^{(n)})^m$. У цьому випадку оператор \eqref{L-Cauchy} сюр'єктивний, а його ядро складається з усіх вектор-функцій $y\in(C^{(n)})^m$ таких, що $(L,C_{m,m})y\in\mathcal{O}\times\mathbb{C}^{m-l}$, де $\mathcal{O}$ позначає нульовий підпростір простору $(C^{(n-1)})^m\times\mathbb{C}^{l}$. Тому на підставі ізоморфізму \eqref{(L,C)isomorphism} робимо висновок, що  вимірність цього ядра дорівнює вимірності простору $\mathcal{O}\times\mathbb{C}^{m-l}$, тобто $m-l$. Отже, оператор  \eqref{L-Cauchy} разом із $(L,B)$ є фредгольмовим з індексом $(m-l)-0$ і в останньому випадку.

Теорему \ref{th_fredh-bis} доведено.

\textit{Доведення теореми $\ref{th_FR_n}$}. Ядро $\ker L$ оператора $L:(C^{(n)})^m\to(C^{(n-1)})^m$ складається з усіх вектор-функцій вигляду $Yq$, де $q\in\mathbb{C}^{m}$. Нагадаємо, що $Y$~--- єдиний розв'язок матричної задачі Коші \eqref{Cachy}. Відображення $C:y\mapsto y(a)$ встановлює ізоморфізм на парі просторів $\ker L$ і $\mathbb{C}^{m}$, причому $C(Yq)=q$ для кожного $q\in\mathbb{C}^{m}$. З означення характеристичної матриці $M(L,B)$ безпосередньо випливає рівність $B(Yq)=M(L,B)q$ для будь-якого $q\in\mathbb{C}^{m}$. Якщо $y\in\ker(L,B)$, то $y=Yq$ для деякого вектора $q\in\mathbb{C}^{m}$ і $By=0$ й тому $M(L,B)q=B(Yq)=0$, тобто $Cy=q\in\ker M(L,B)$. Зворотно, якщо $q\in\ker M(L,B)$, то $B(Yq)=M(L,B)q=0$ й тому $y:=Yq\in\ker(L,B)$, причому $Cy=q$. Таким чином, ізоморфізм $C:\ker L\leftrightarrow\mathbb{C}^{m}$ встановлює взаємно однозначну відповідність між $\ker(L,B)$ і $\ker M(L,B)$. Це дає потрібну формулу \eqref{dimker}, тобто
\begin{equation*}
\dim\ker(L,B)=\dim\ker M(L,B)=m-\operatorname{rank}M(L,B).
\end{equation*}
З неї випливає формула \eqref{dimcoker1} на підставі теореми~\ref{th_fredh-bis}. Теорему~\ref{th_FR_n} доведено.

\textit{Доведення теореми $\ref{th_hm_stc}$}.
Припустимо, що $(L_k,B_k)\xrightarrow{s}(L, B)$; усі границі розглядаємо при $k\to\infty$. Тоді $A_{k}z_{j}=L_{k}z_{j}\to Lz_{j}=A z_{j}$ у просторі $(C^{(n-1)})^{m}$ для кожного $j\in\{1,\ldots,m\}$, де $z_{j}(t):=(\delta_{j,1},\ldots,\delta_{j,m})$ при $a\leq t\leq b$; як звичайно, $\delta_{j,i}$~--- символ Кронекера. Отже, $j$-ий стовпчик матриці-функції $A_{k}$ прямує до $j$-го стовпчика матриці-функції $A$ у цьому просторі, тобто $A_{k}\to A$ в $(C^{(n-1)})^{m\times m}$. 

Розглянемо лінійні відображення $\mathcal{L}:Z\mapsto Z'+AZ$ і $\mathcal{L}_{k}:Z\mapsto Z'+A_{k}Z$, де $Z\in(C^{(n)})^{m\times m}$. Вони є обмеженими операторами на парі просторів 
\begin{equation}\label{L-operators}
\mathcal{L},\mathcal{L}_{k}:
(C^{(n)})^{m\times m}\to(C^{(n-1)})^{m\times m}.
\end{equation}
Для довільної матриці-функції $Z\in(C^{(n)})^{m\times m}$ маємо:
\begin{gather*}
\|\mathcal{L}_{k}Z-\mathcal{L}Z\|_{(n-1)}=\|(A_{k}-A)Z\|_{(n-1)}\leq\\
\leq c\,\|A_{k}-A\|_{(n-1)}\|Z\|_{(n-1)}
\leq c\,\|A_{k}-A\|_{(n-1)}\|Z\|_{(n)},
\end{gather*}
де число $c>0$ залежить лише від $n$ і $m$ (бо простір $(C^{(n-1)})^{m\times m}$ є банаховою алгеброю). Тому
\begin{equation*}
\frac{\|\mathcal{L}_{k}Z-\mathcal{L}Z\|_{(n-1)}}{\|Z\|_{(n)}}\to 0,
\end{equation*}
тобто для операторів \eqref{L-operators} виконується збіжність $\mathcal{L}_{k}\rightrightarrows\mathcal{L}$ у рівномірній операторній топології.  

Розглянувши лінійне відображення $\mathcal{C}:Z\mapsto Z(a)$ на вектор-функціях $Z\in(C^{(n)})^{m\times m}$, отримаємо топологічні ізоморфізми
\begin{equation*}
(\mathcal{L},\mathcal{C}),(\mathcal{L}_{k},\mathcal{C}):
(C^{(n)})^{m\times m}\leftrightarrow
(C^{(n-1)})^{m\times m}\times\mathbb{C}^{m}
\end{equation*}
згідно з теоремою про існування і єдиність розв'язку задачі Коші. Оскільки  $(\mathcal{L}_{k},\mathcal{C})\rightrightarrows(\mathcal{L},\mathcal{C})$ за доведеним, то для обернених операторів виконується збіжність  $(\mathcal{L}_{k},\mathcal{C})^{-1}\rightrightarrows
(\mathcal{L},\mathcal{C})^{-1}$. Тому
\begin{equation*}
Y_{k}:=(\mathcal{L}_{k},\mathcal{C})^{-1}(O_{m},I_{m})\to (\mathcal{L},\mathcal{C})^{-1}(O_{m},I_{m})=Y
\end{equation*}
у просторі $(C^{(n)})^{m\times m}$. Тут $Y$ і $Y_{k}$~--- єдині розв'язки відповідно матричних задач Коші \eqref{Cachy} і 
\begin{equation*}
Y_{k}'(t)+A_{k}(t)Y_{k}(t)=O_m, \quad t \in (a,b), \quad\quad Y_{k}(a)=I_m.
\end{equation*}
Оскільки $B_k\xrightarrow{s}B$ за умовою, то із встановленої збіжності $Y_{k}\to Y$ в $(C^{(n)})^{m\times m}$ випливає за означенням характеристичної матриці, що $M(L_k,B_k)\to M(L,B)$.

Теорему~$\ref{th_hm_stc}$ доведено.

\textit{Доведення теореми~$\ref{th_stc-col}$}.
Для числової матриці $M(L,B)$ існує ненульовий мінор порядку $\operatorname{rank}M(L,B)$. Оскільки $M(L_k,B_k)\to M(L,B)$ при $k\to\infty$ за теоремою~$\ref{th_hm_stc}$, то відповідний мінор матриці 
$M(L_k,B_k)$ відмінний від нуля для всіх $k\gg1$. Отже, $\operatorname{rank}M(L_k,B_k)\geq\operatorname{rank}M(L,B)$ для всіх $k\gg1$. Звідси на підставі теореми~\ref{th_FR_n} випливають потрібні формули \eqref{inq ker} і \eqref{inq coker}. Теорему~$\ref{th_stc-col}$ доведено.

\section{Acknowledgements}\label{section6}

Робота виконана у рамках науково-дослідної теми молодих учених НАН України 2024--2025 рр. (державний реєстраційний номер~--- 0124U002111). Автор також підтриманий стипендією Президента України для молодих учених, грантом Simons Foundation (1290607, V.O.S.) і грантом імені Марії Склодовської-Кюрі №~873071 за програмою Європейського Союзу Горизонт 2020~--- Рамкова програма з досліджень та інновацій: Спектральна оптимізація: від математики до фізики і передових технологій.

\bigskip

\textbf{Vitalii Soldatov}

{Tereshchenkivska Str. 3, 01024 Kyiv, Ukraine}

{Institute of Mathematics of NAS of Ukraine, Kyiv 
	
Institute of Applied Mathematics and Mechanics of NAS of Ukraine, Sloviansk
}

\textbf{В. О. Солдатов }

01024, Київ, вул.~Терещенківська, 3

Інститут математики НАН України, Київ

Інститут прикладної математики і механіки НАН України, Слов'янськ

{soldatov@imath.kiev.ua, soldatovvo@ukr.net}

\end{document}